# Concordance probability in a big data setting: application in non-life insurance

Robin Van Oirbeek
*Data Office*
*Allianz Benelux*
Brussels, Belgium
robin.vanoirbeek@gmail.com

Christopher Grumiau
*Data Office*
*Allianz Benelux*
Brussels, Belgium
christopher.grumiau@allianz.be

Tim Verdonck
*Department of Mathematics.*
*University of Antwerp*
Antwerp, Belgium
tim.verdonck@uantwerp.be

*Abstract*— The concordance probability or C-index is a popular measure to capture the discriminatory ability of a regression model. In this article, the definition of this measure is adapted to the specific needs of the frequency and severity model, typically used during the technical pricing of a non-life insurance product. Due to the typical large sample size of the frequency data in particular, two different adaptations of the estimation procedure of the concordance probability are presented. Note that the latter procedures can be applied to all different versions of the concordance probability.

## I. Introduction

When determining the premium of a non-life insurance product, such as a car insurance product, its technical tariff is mostly used as a starting point. This technical tariff captures how the insured risk depends on the customer profile, which can be gauged by covariates such as the customer's age, his/her driving experience, the car type he/she is driving, and so on. Typically, a Poisson regression model is used to model the frequency or number of accidents that a customer will have during the considered period (also defined as the exposure), and a Gamma regression model for the severity of the given accident.

Constructing such a technical tariff often is a non-linear and labor-intensive process. One frequently needs to decide which covariates should be retained in the model, which interaction effects should be included, etc. Significance of the considered effect or covariate generally determines whether or not the effect or covariate will be retained. We might want to step away from just using significance as a variable selection criterion by using the concordance probability or C-index as a selection criterion. Indeed, distinguishing risky clients from less risky clients is at the core of what insurance is all about and moreover significance does not necessarily capture the model's ability to do so.

The concordance probability measures the discriminatory ability of a model, i.e. the consistency of the model predictions and observed outcomes [1]. However, the basic version of this measure is not equipped to deal with the concept of exposure, i.e. the duration of a policy or insurance contract, which plays a pivotal role in frequency models. Therefore, we will generalize the concordance probability to cope with the concept of exposure as well. Moreover, the extension of the concordance probability [2] will be adapted to the severity model.

Data sets modelling the frequency are typically rather large, mostly comprising of multiple millions of observations. Since estimating the concordance probability is computationally intensive, it would take hours or even days to calculate it on such large data sets, in a naive manner.

In this paper, we propose two approximations to efficiently compute the concordance probability for big data. The first one takes a sample of the data set, for which the concordance probability is calculated, and computes a confidence interval to gauge the precision of the obtained estimate. The second solution makes use of the $k$-means clustering algorithm to split customers up in different groups, which significantly accelerates the computation of the concordance probability.

In Section II, the basic definition of the concordance probability is presented. The adaptations of this basic definition to the frequency and severity model, as well as their estimation in a big data setting, are presented in Section III. The newly developed measures are applied to a real data set in Section IV, followed by some conclusions in Section V.

## II. Literature Review

The concordance probability is commonly used to assess the predictive ability of a binary regression model and equals the area under the ROC curve [1]. More specifically, it corresponds to the probability that a randomly selected subject with outcome $Y = 1$ has a higher predicted probability $\pi(X) = P(Y = 1 | X)$ than a randomly selected subject with outcome $Y = 0$, where $X$ corresponds to the vector of covariates [3]. This is equivalent to the probability that a randomly selected comparable pair $(i, j)$, i.e. a random pair with $Y_i = 1$ and $Y_j = 0$, is a concordant pair, i.e. a pair with $\pi(X_i) > \pi(X_j)$:

$$C = P(\pi(X_i) > \pi(X_j) | Y_i = 1 \& Y_j = 0). \quad (1)$$

The concordance probability $C$ normally ranges between 0.5 and 1, and the closer it is to 1, the better its discriminatory ability. If its value drops below 0.5, the predictions are consistently inconsistent. For a sample of size $n$, the concordance probability typically is estimated as:

$$\hat{C} = \frac{n_c}{n_t} = \frac{\hat{\pi}_c}{\hat{\pi}_c + \hat{\pi}_d}$$

$$= \frac{\sum_{i=1}^{n-1} \sum_{j=i+1}^{n} I(\hat{\pi}(X_i) > \hat{\pi}(X_j) \& y_i = 1 \& y_j = 0)}{\sum_{i=1}^{n-1} \sum_{j=i+1}^{n} I(y_i = 1 \& y_j = 0)}, \quad (2)$$

corresponding to the ratio of the number of concordant pairs $n_c$ over the number of comparable pairs $n_t$. The value $\hat{\pi}_c$ ($\hat{\pi}_d$) refers to the estimated probability that a comparable pair is concordant (discordant) respectively and $I(.)$ to the indicator variable. Note that $y_i$ and $\hat{\pi}(X_i)$ correspond to the observed outcome and estimated predicted probability for observation $i$ with $X_i$ the vector of observed covariates respectively. For more information on the different aspects of the predictive ability, such as the difference between the discriminatory ability and the calibration, we refer to [4].

We assume that ties on the observed and predicted values preferably are excluded from the calculation of any concordance probability as ties attenuate the concordance probability to 0.5, depending on the distribution of the observations over the different risk groups, which is undesirable [5]. See [6] for a more detailed treatment of the subject of ties on the concordance probability.

## III. METHODOLOGY

The basic definition, presented in the previous section, will now be extended to the specific needs of the data used for the technical pricing exercise. In this generalization, the total cost for a non-life insurance contract of a given duration or exposure is estimated for a given customer profile. Two data sets are necessary to this end: frequency data (i.e. the number of events, the exposure and characteristics of the policyholder), and severity data (i.e. the cost of an event and characteristics of the policyholder).

### A. Generalizations for frequency data

The basic definition of $C$ requires the definition of 2 groups, based on the number of events that occurred during the duration of the policy. Non-life insurance contracts, such as car and home insurance contracts, typically have an exposure of 1 year or less, such that more than 2 events infrequently take place during this (short) period. Hence, 3 groups will be defined: policies that experienced 0 events, 1 event, and 2 events or more, represented by the 0, 1 and 2 group, resulting in the following 3 definitions of the concordance probability:

$$C_{01+} = P(\pi^N(X_i) < \pi^N(X_j) \mid Y_i^N = 0 \ \& \ Y_j^N \geq 1),$$

$$C_{02+} = P(\pi^N(X_i) < \pi^N(X_j) \mid Y_i^N = 0 \ \& \ Y_j^N \geq 2), \quad (3)$$

$$C_{12+} = P(\pi^N(X_i) < \pi^N(X_j) \mid Y_i^N = 1 \ \& \ Y_j^N \geq 2),$$

where $\pi^N(.)$ refers to the predicted frequency of the chosen frequency model, typically being a Poisson model, and $Y^N$ corresponds to observed claim number. Definitions (3) have interesting interpretations: $C_{01+}$ ($C_{02+}$) evaluates the ability of the model to discriminate policies that did not experience accidents from policies that experience at least 1 (2) accident(s) respectively. $C_{12+}$ quantifies the ability of the model to discriminate policies that experienced 1 accident from policies that experienced multiple accidents, hence of clients that could just have been unfortunate versus clients that are (probably) accident-prone.

In order to make sure that the pair is comparable, the exposure should be taken into account in the above definitions. Indeed, the exposure $\lambda$ of both members of the pair should be more or less the same to be fairly comparable. Since $\lambda$ is (technically) a number of the interval $[0,1] \subseteq \mathbb{R}$, it does not make sense to impose that the exposure of both elements of the pair is exactly the same. The definitions of (2) are therefore equal to:

$$C_{01+}^{\approx} = P(\pi^N(X_i) < \pi^N(X_j) \mid Y_i^N = 0 \ \& \ Y_j^N \geq 1 \ \& \ \lambda_i \approx \lambda_j),$$

$$C_{02+}^{\approx} = P(\pi^N(X_i) < \pi^N(X_j) \mid Y_i^N = 0 \ \& \ Y_j^N \geq 2 \ \& \ \lambda_i \approx \lambda_j), (4)$$

$$C_{12+}^{\approx} = P(\pi^N(X_i) < \pi^N(X_j) \mid Y_i^N = 1 \ \& \ Y_j^N \geq 2 \ \& \ \lambda_i \approx \lambda_j).$$

In practice, the maximal difference in exposure between both members of a pair that is considered to be negligible, is a tuning parameter. A local version of the above definitions can be introduced as well:

$$C_{01+}^{\approx}(\lambda) = P(\pi^N(X_i) < \pi^N(X_j) \mid Y_i^N = 0 \ \& \ Y_j^N \geq 1 \ \& \ \lambda_i \approx \lambda_j \approx \lambda)$$

$$C_{02+}^{\approx}(\lambda) = P(\pi^N(X_i) < \pi^N(X_j) \mid Y_i^N = 0 \ \& \ Y_j^N \geq 2 \ \& \ \lambda_i \approx \lambda_j \approx \lambda), (5)$$

$$C_{12+}^{\approx}(\lambda) = P(\pi^N(X_i) < \pi^N(X_j) \mid Y_i^N = 1 \ \& \ Y_j^N \geq 2 \ \& \ \lambda_i \approx \lambda_j \approx \lambda),$$

where $\lambda$ is a tuning parameter, corresponding to the exposure value for which the local concordance probability needs to be computed. The appealing aspect of definitions (5) is that it allows us to construct a $(\lambda, C(\lambda))$ plot, i.e. an evolution of the local concordance probabilities in function of the exposure. After all, the global concordance probability of (4) can be seen as an average value of the local concordance probability.

### B. Generalization for severity data

Since it might be of little practical importance to distinguish claims from one another that only slightly differ in claim cost, the basic definition (1) can be extended as:

$$C(v) = P(\pi^S(X_i) < \pi^S(X_j) \mid Y_i^S < Y_j^S \ \& \ Y_j^S - Y_i^S \geq v), (6)$$

where $\pi^S(.)$ refers to the predicted claim size of the chosen severity model, typically being a Gamma model, and $Y^S$ corresponds to observed claim size. In other words, one just takes those claims into account of which the claim size differs with at least a value $v$ from one another, hereby selecting pair of claims that make more sense from a business point of view. Also, a $(v, C(v))$ plot can be constructed where different values for the threshold $v$ are chosen, as to investigate the influence of $v$ on (6). Note that, as shown in [1], extra conditions can be added to definition (6).

### C. Point Estimate and Confidence Interval

Each of the presented concordance probabilities can be estimated in a similar style as shown in (2), i.e. by taking the ratio of the number of concordant pairs over the number of comparable pairs. When dealing with even a moderate sized data set, it can take a considerable time to estimate the concordance probability. This is especially so for frequency data since its sample size easily exceeds 100,000 observations, even for small insurance companies. Therefore strategies need to be elaborated to deal with this big data setting. Probably, the easiest manner is by sampling at random observations from the full data set and to check against the remaining observations of the full data set how many concordant and comparable pairs can be enumerated. The algorithm hence corresponds to:

1. Choose a value for $S$, corresponding to the number of observations that will be randomly sampled (without replacement) from the full data set.

2. Select randomly one observation $(i)$ from the full data set.

3. Determine the number of concordant and comparable pairs, regarding the remaining observations in the full data set and the selected observation $i$, resulting in $n_{c,i}^*$ and $n_{t,i}^*$ respectively.
4. Remove observation $i$ from the full data set and repeat the algorithm again starting from step 2 as long as the number of selected observations is lower than $S$.

Given that $S$ generally is much smaller than $n$ and yet still sufficiently large, no involuntary bias should be caused by the above estimation method. Hence,

$$\hat{C} = \frac{\hat{\pi}_c}{\hat{\pi}_c + \hat{\pi}_d} \approx \frac{\sum_{i=1}^{S} n_{c,i}^*}{\sum_{i=1}^{S} n_{t,i}^*}. \quad (7)$$

An alternative method would be by taking a sample of size $S$ of the full data set, and to just evaluate the number of comparable and concordant pairs for the observations of the latter sample only. Despite the computing time being still within reasonable bounds on modern-day computers, i.e. a couple of minutes, the precision of the resulting estimates was less high than the one obtained from the above procedure.

In general, for the above procedure, good values for $S$ lie within the range $[5,000; 20,000]$, as informally tested on several data sets using modern-day computers. In any case, if one chooses $S$ too low, meaning that the approximation is too coarse, a higher value should be chosen and the resulting precision should be checked. The latter procedure can be repeated until a satisfactory precision is obtained. This precision can be gauged by constructing a confidence interval for the latter estimate. As such, by tuning sample size $S$, the user can determine her/himself when satisfactory precision is obtained.

In [3], a general methodology to obtain the confidence interval for any definition of the concordance probability is presented:

$$\hat{C} \pm z_{\alpha/2} \sqrt{\frac{\widehat{var}[\hat{C}]}{n}},$$

$$\widehat{var}[\hat{C}] = \frac{4(\hat{\pi}_d^2 \hat{\pi}_{cc} - 2\hat{\pi}_c \hat{\pi}_d \hat{\pi}_{cd} + \hat{\pi}_c^2 \hat{\pi}_{dd})}{(\hat{\pi}_c + \hat{\pi}_d)},$$

where $\alpha$ equals the confidence level, typically equal to 0.05, and $\hat{\pi}_{cc}$ ($\hat{\pi}_{dd}$) corresponds to the estimated probability that a randomly selected observation is concordant (discordant) with 2 randomly selected observations. $\hat{\pi}_{cd}$ is the estimated probability that a randomly selected observation is concordant with 1 randomly selected observation and discordant with 1 other randomly selected observation. Note that [3] assumes that each pair is ordered and ties are included, which explains the slight adaptations of our estimators:

$$\hat{\pi}_{cc} = \sum_{i=1}^{S} \left(\frac{n_{c,i}^*}{n_{t,i}^*}\right)\left(\frac{n_{c,i}^* - 1}{n_{t,i}^* - 1}\right)\left(\frac{n_{t,i}^*}{n_{t,S}^*}\right) \approx \frac{1}{n_{t,S}^*} \sum_{i=1}^{S} \frac{(n_{c,i}^*)^2}{n_{t,i}^*},$$

$$\hat{\pi}_{dd} \approx \frac{1}{n_{t,S}^*} \sum_{i=1}^{S} \frac{(n_{t,i}^* - n_{c,i}^*)^2}{n_{t,i}^*},$$

$$\hat{\pi}_{cd} \approx \frac{1}{n_{t,S}^*} \sum_{i=1}^{S} \frac{(n_{t,i}^* - n_{c,i}^*) n_{c,i}^*}{n_{t,i}^*},$$

where $n_{t,S}^* = \sum_{i=1}^{S} n_{t,i}^*$. Hence, if the width of the confidence interval is too large, $S$ should be increased, and the concordance probability should be recalculated using (7) until a satisfactory width is obtained. As a rule of thumb, a width of 0.01 to 0.02 is generally acceptable.

*D. k-means approximation*

As can be seen in (1) and (3), each definition of a concordance probability always consists of 2 groups. These are defined based on the observed outcomes, of which the distributions of their respective predictions are tested on consistency. As a result, by approximating the distribution of the predictions of both groups in one way or the other, a faster way for computing an approximation of the concordance probability can be obtained. In the remainder of this section, 2 other approximations will be presented.

The $k$-means approximation can be obtained by applying a clustering algorithm to the predictions of each group in the definition of the concordance probability separately. The $k$-means algorithm is a good candidate as it scales linearly with the number of observations, such that it is well suited to be used in a big data setting. The number of clusters $k$ will then determine the level of approximation, meaning that a more precise estimate will be obtained as $k$ increases. Note that the latter methodology could easily be replaced by another clustering algorithm besides the $k$-means clustering algorithm.

After having applied the clustering algorithm, only the cluster centroids will be considered in estimating the concordance probability. Indeed, all clusters of each group will be compared with one another and the importance of each cluster by cluster comparison is weighted by the probability that a randomly selected pair of observations belongs to the respective clusters, or:

$$\hat{C} = \frac{\hat{\pi}_c}{\hat{\pi}_c + \hat{\pi}_d} \approx \sum_{i=1}^{k-1} \sum_{j=i+1}^{k} I\left(\widehat{\pi_B^i} > \widehat{\pi_A^j}\right) w_B^i w_A^j \quad (8)$$

where subscript B (A) refers to the group of observations of which the predictions are supposed to be larger (lower) than the ones of subscript A (B) respectively. $\widehat{\pi_*^l}$ is the centroid of the $l$-th cluster of group $*$ and $w_*^l$ the weight of the $l$-th cluster of group $*$. The latter is estimated by computing the observed percentage of observations of the group $*$ that belong to cluster $l$, such that by definition $\sum_{l=1}^{k} w_*^l = 1$.

When applying the $k$-means approximation (8) to definitions (4) and (5), clusters will only be determined for observations that pertain to a certain range of exposure. The above methodology will then be applied unaltered to this subset of the whole data set only and repeated for each subset separately. The obtained estimates will then be combined in a similar manner as presented by definitions (4) and (5).

IV. FINDINGS

In this section, the presented measures and approximations will be tested on a publicly available insurance data set, consisting of 159,911 observations, 3 continuous covariates and 4 categorical covariates. 145,719

(or 91.1%) policyholders did not file a claim, 12,874 (or 8.05%) filed only 1 claim and 1,354 (or 0.85%) filed 2 or more claims. For the frequency analysis, a Poisson regression model was used, whereas for the claim size analysis a Gamma regression model was used.

Due to the large sample size of the considered data set, the estimation method testing all possible comparable pairs (2) was not tested out, only estimation methods (7), (8) and (10) were considered. We start by using estimation method (7) to explore the predictive ability of both Poisson and Gamma models in computing definitions (4), (5) and (6).

In Fig. 3, the evolution of the predictive ability in function of the exposure is shown for the Poisson model with main effects only, i.e. definition (5). A maximum absolute difference of 0.05 in exposure between both members of the comparable pairs was chosen for definitions (4) and (5). Finally, $S$ was set to 20000. The computational time was around 220 seconds for $C^{\approx}_{01+}$ and $C^{\approx}_{02+}$ and around 125 seconds for $C^{\approx}_{12+}$. We can see that the estimates of all concordance probabilities remain quite stable in function of the exposure, except for $C^{\approx}_{12+}(\tau)$ which is rather low for small exposures, but which increases as the exposure increases. The global versions of the 3 considered concordance probabilities are quite similar to values for the local version at exposure value 1, or $C^{\approx}_{ij+}(1) \approx C^{\approx}_{ij+}$ with $i$ and $j \in \{0,1,2\}$. This is due to 27.81% of the observations having an exposure of one. Note that $C^{\approx}_{12+} < C^{\approx}_{01+} \ll C^{\approx}_{02+}$, and that it is hard for the model to distinguish policy holders who have had just one accident from policy holders who have had multiple accidents ($C^{\approx}_{12+} = 0.593$), and that the model distinguishes a bit better policy holders who have not had an accident from those who have had at least one accident ($C^{\approx}_{01+} = 0.621$). However, the model is quite good at separating policy holders who have not had an accident from policy holders who have had many accidents ($C^{\approx}_{02+} = 0.690$). This is a pattern that was observed by the authors on other non-life insurance data sets as well.

The impact of modelling the covariate effects of the 3 continuous variables using splines is shown in Fig. 4: the global and the local concordance probabilities are a bit higher than what can be observed in Fig. 3. The impact is rather limited ($\Delta C^{\approx}_{01+} = 0.003, \Delta C^{\approx}_{12+} = 0.001, \Delta C^{\approx}_{02+} = 0.013$), but it is non negligible nevertheless. These improvements can therefore almost completely be attributed to model misspecification in the main effects model. The evolution of the predictive ability in function of the exposure is quite similar for both models, yet the curves are a bit higher in the lower exposure area for the spline model.

For the considered global estimates of the concordance probability, the confidence interval (CI) is computed to check whether or not setting $S$ to 20000 was sufficiently large or not. In Tab. 1, one can see that the differences between the lower and upper bound of the CI is around 0.015, which is acceptable in practice.

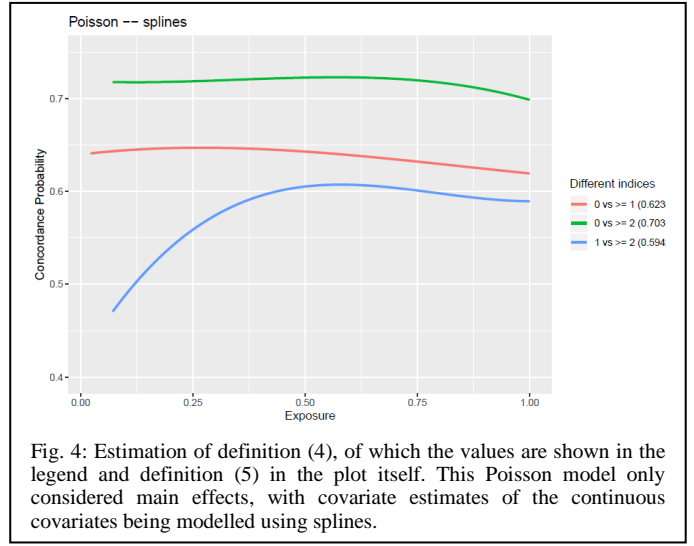

Fig. 4: Estimation of definition (4), of which the values are shown in the legend and definition (5) in the plot itself. This Poisson model only considered main effects, with covariate estimates of the continuous covariates being modelled using splines.

The predictive ability of the Gamma model is investigated next. The latter model consists of main effects of the 7 covariates, with no spline estimates for the covariate effects of the continuous covariates. The evolution of $C(v)$ is shown in Fig. 5 and one clearly sees that the discriminatory ability is quite bad when $v$ is low, meaning that practically all possible pairs are allowed. As $v$ increases, the pairs where there is little difference in observed claim size are removed, such that less similar claims remain. This results in an uplift of $C(v)$ as $v$ increases. The pointwise CI shows that $S$ equals to 5000 suffices to compute $C(v)$ values with an acceptable precision.

The last analysis evaluates the performance of the $k$-means approximation on the data set. There are 2 tuning parameters, of which is the first one is the number of splits in the exposure distribution, necessary to comply the methodology with definition (4). The other tuning parameter is the number of clusters $k$. Note that the $k$-means approximation necessitates to be rerun with other initial values a certain number of times when the number of clusters and/or the number of splits of the exposure distribution is insufficient. This is due to the non-negligible variance over the estimator when these reruns are not performed (see Tab. 2).

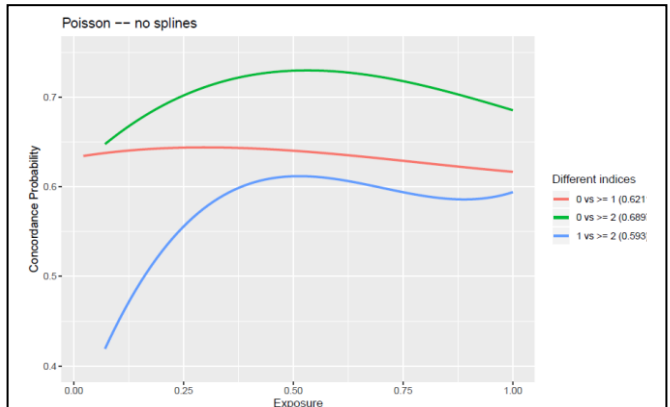

Fig 3: Estimation of definition (4), of which the values are shown in the legend and definition (5) in the plot itself. This Poisson model only considered main effects, with the covariate estimates of the continuous covariates not being modelled using splines.

Tab. 1. The estimates together with the corresponding CI of the global concordance probabilities $C^{\approx}_{01+}, C^{\approx}_{02+}$ and $C^{\approx}_{12+}$ for the Poisson model with and without splines on the continuous covariate effects.

|  |  | *Estimate* | *Lower Bound CI* | *Upper Bound CI* |
|---|---|---|---|---|
| Without splines | $C^{\approx}_{01+}$ | 0.621 | 0.614 | 0.629 |
|  | $C^{\approx}_{02+}$ | 0.690 | 0.682 | 0.697 |
|  | $C^{\approx}_{12+}$ | 0.593 | 0.585 | 0.601 |
| With splines | $C^{\approx}_{01+}$ | 0.624 | 0.616 | 0.631 |
|  | $C^{\approx}_{02+}$ | 0.703 | 0.695 | 0.710 |
|  | $C^{\approx}_{12+}$ | 0.594 | 0.582 | 0.597 |

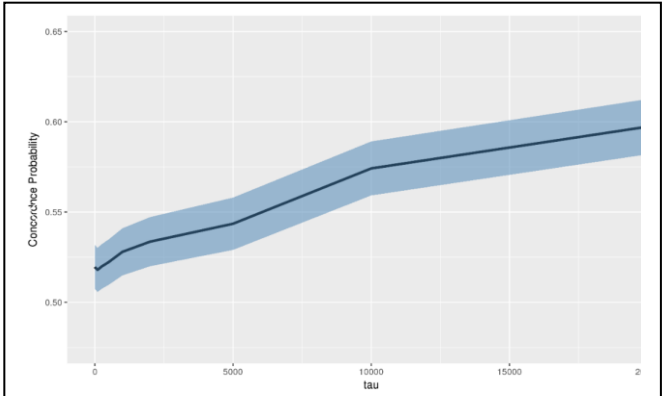

Fig. 5: Estimation of definition (6), in function of $v$. A pointwise CI was added to the plot to assess the precision of the resulting estimates.

From Tab. 2, one can see that relatively good to excellent estimates were obtained for all investigated scenarios. Some patterns could be discerned: the $k$-means approximation tend to underestimate the concordance probability when the approximation is too coarse. Note that sometimes not that many reruns are necessary, and this especially for the finer approximations, hence the ones that require a much longer computational time per rerun. However, to be on the safe side, one could always opt for a large number of reruns. In that case the computational time will strongly increase, especially for the finer approximations

## V. Conclusion

In this article, the basic definition of the concordance probability is extended to the special needs of the technical pricing of a non-life insurance product. Since data sets in insurance are generally rather large, the two novel estimation methods are proposed to compute the concordance probability in a big data context, all resulting in approximations of the desired concordance probability. The performance of these approximations is evaluated on a publicly available insurance data set. In general, both estimators resulted in good estimation properties, however, the $k$-means approximation was computationally significantly less intensive than the first approximation, and is therefore recommended by the authors. In future research, an extensive simulation study should be set up to investigate the estimation properties of both estimation methods in more detail. The estimation properties of the proposed CI could be investigated in more detail, and a CI could be constructed and investigated for the $k$-means approximation.

Tab. 2. Table showing the estimated concordance probability $C_{01+}^{\approx}$ for the Poisson model without spline estimates on the continuous covariate effects and within brackets the computational time in seconds for the $k$-means approximation, in function of the exposure splits (8, 15 and 70) and the number of clusters (10, 19 or 50). As shown in the second column for the $k$-means approximation, the approximation with 10 and 19 clusters necessitated 50 and 20 reruns respectively to obtain stable estimates. The best solution is shown in bold.

| $C_{01+}^{\approx}$ | | 8 | 15 | 70 |
|---|---|---|---|---|
| $k$-means | 10 (50) | 0.615 (9.8) | 0.608 (11.5) | 0.603 (27.1) |
| | 19 (20) | 0.626 (3.4) | 0.619 (6.9) | 0.620 (15.8) |
| | 50 | 0.626 (1.2) | 0.621 (1.9) | **0.619 (2.2)** |